\documentclass{article}

\usepackage{amssymb,amsmath} %%%

\usepackage{amscd}%%commutative diagrams

\usepackage[latin1]{inputenc} \usepackage{times}
\usepackage{amssymb,amsmath} \usepackage{latexsym} \usepackage{amsthm}

\newtheorem{theorem}{Theorem}[section]

\newtheorem{proposition}[theorem]{Proposition}

\newfont{\Goth}{eufm10 at 11pt} \newfont{\biggoth}{eufm10 at 18pt}

\newcommand{\HH}{{\mathcal H}}

 \newcommand{\mgoth}[1]{\mbox{\Goth #1}}
 \newcommand{\C}{\mathbb C}
 \newcommand{\R}{\mathbb R}

\newcommand{\gra}{\alpha}

\title{\bf  K-contact Lie groups of dimension five or greater}
\author{Brendan J. Foreman}

\begin{document}

\maketitle

\begin{abstract}{ \footnotesize
We prove that a K-contact Lie group of dimension five or greater is the central extension of a symplectic Lie group by complexifying the Lie algebra and applying a result from complex contact geometry, namely, that, if the adjoint action of the complex Reeb vector field on a complex contact Lie algebra is diagonalizable, then it is trivial.

}
\end{abstract}

\section{Introduction}

Recall that a {\it real contact structure} on a manifold $M$ of dimension $2n+1$ is a distribution $\HH$ of $TM$ given as the kernel of a 1-form $\eta$ satisfying $\eta\wedge d\eta^n\neq 0$ at all points of $M$.  The {\it Reeb vector field} of a contact manifold
$(M, \HH, \eta)$ is the vector field $\xi$ transverse to $\HH$ defined by the equations
$$\eta(\xi)=1,\ \iota(\xi)dn=0.$$ The tangent bundle of $M$ splits by $TM=\HH \oplus \left<\xi\right>,$ and we denote the projection $TM \rightarrow \HH$ by $\HH,$ as well.
If $M$ is a Lie group such that $\eta$ is left-invariant, then we call $M$ a contact Lie group.  

Note that all of the above definitions also make sense if we switch to the complex category.  That is, we call $G$ a {\it complex contact Lie group}, if $G$ is a complex Lie group with a left-invariant holomorphic one-form $\eta$ such that $\eta\wedge\eta^n\neq0$ for $dim_\C G=2n+1.$  Similarly, the definitions of the complex contact distribution and Reeb vector field carry over analogously.  

The main result of this paper, namely that a K-contact Lie group of dimension five or greater is the central extension of a sympletic Lie group is the result of this analogy.  Namely, given the K-contact Lie group, we complexify the contact structure, use a result in complex contact geometry and then note the consequences on the original real contact Lie group.  Interestingly, this is the same strategy for which twistor spaces were originally invented and utilized in \cite{penrose}. See \cite{blair} for additional and more detailed information on both real and complex contact structures.

\section{Real contact metric structures}

This section provides the preliminary definitions and results with the real contact geometry.  A metric $g$ on a contact manifold $(M,\eta)$ is
called {\it associated} if the following criteria are satisified
\begin{enumerate} \item $\eta(X)=g(X, \xi)$ for all $X\in TM$ and \item
the endomorphism $\phi:TM\rightarrow TM$ defined for $X,\ Y\in TM$ by
$$g(X,\phi Y)=d\eta(X,Y)$$ satisfies $$\phi^2=-I+\eta\otimes\xi.$$
\end{enumerate} 
Much is known about the resulting Riemannian geometry of associated metrics on contact manifolds (\cite{blair}).  For the purposes here, the following results are needed.

\begin{proposition}Let $(M, \HH, \eta)$ be a contact manifold of dimension $2n+1$ with associated metric $g.$  Then the Levi-Civita connection $\nabla$ satisfies
$$\nabla_X \xi=-\phi X - \phi hX,$$
where $\phi$ is a skew-symmetric endomorphism of $TM$ such that $\phi^2=-Id+\eta\otimes\xi$ and $h$ is symmetric with respect to $g.$
\end{proposition}

\noindent For symplectic manifolds,
there is an analogous concept of associated metric, namely, a metric $k$ is {\it associated} to the symplectic structure of a manifold $S$, if there is an almost complex structure $J$ on $S$ such that the symplectic form $\omega$ is given by $\omega(X,Y)=g(X,JY).$  

An associated metric $g$ of a contact manifold $(M, \HH, \eta)$
is called {\it K-contact}, if $\xi$ is an infinitesimal automorphism of
$g,$ i.e., $\mathcal L_\xi g=0.$ It is not difficult to see that this is equivalent to the nullity of the tangent bundle transformation $h$ as given in the proposition above.  Also, it is easy to see that, if there is a symplectic manifold
$(S,\omega)$ such that $\pi: M\rightarrow S$ is a fibration of the
leaves of the Reeb vector field with $\pi^* \omega= d\eta,$ then an associated metric $g$ on $M$ is
K-contact if and only if there is an associated metric $k$ on
$(S,\omega)$ such that $\pi^*(k)_{|\HH}=g_{|\HH}.$ (see \cite{blair}).

\begin{proposition} Let $G$ be a contact Lie group with left-invariant
contact form $\eta$, Reeb vector field $\xi$ and left-invariant
associated metric $g.$  Then $g$ is K-contact if and only the matrix
form of $ad(\xi)$ on the Lie algebra \mgoth{g}  of $G$ is skew-symmetric
with respect to any orthonormal basis $e=\{e_1,\dots e_{2n}\}$ of the
contact distribution $\HH=ker\eta.$ \end{proposition}

\noindent{\bf Proof:}\ \ Let $g$ be a left-invariant metric and $X,\ Y\,\ Z$ be left-invariant vector fields on $G.$   Then $$g(\nabla_XZ,Y)=-\frac{1}{
2}\left(g([Z,Y],X)+g([X,Y],Z)+g([Z,X],Y)  \right)$$ so that
$$g(\nabla_XZ,Y)+g(\nabla_YZ,X)=-g([Z,Y],X) - g([Z,X],Y) .$$ If $g$ is
associated, then $\nabla_X\xi=-\phi X-\phi hX$ and $g$ is K-contact if
and only if $h=0.$  But the transformation $\phi h$ is the symmetric
part of $X\mapsto \nabla_X \xi.$  So, $h=0$ if and only if
$0=-g([\xi,Y],X) - g([\xi,X],Y) $ for any left-invariant horizontal
vector fields $X$ and $Y,$ i.e., $0=-g(ad(\xi)Y,X) - g(ad(\xi)X,Y) $ for
any $X,\ Y\in \mgoth{g}.$  This proves the proposition.

\vskip .25in It is well known that any real skew-symmetric $n\times n$
matrix $B$ is diagonalizable in the space of complex matrices,
$M_{n\times n}(\C).$  More specifically, there is a $Q\in O(n)$ such
that $$QBQ^t=\left(\begin{array}{c|c|c|c} \begin{array}{cc}0 & b_1\\
-b_1&0\\ \end{array}&&&\\ \hline &\ddots&&\\ \hline &&\begin{array}{cc}0
& b_k\\ -b_k&0\\ \end{array}\\ \hline  &&&0\\ \end{array}\right),$$
\noindent for some $b_1, \dots, b_k\in \R^*.$  Thus, if the Jordan
canonical form of $ad(\xi)$ with respect to any left-invariant basis of
{\Goth g} contains a block matrix of the form
$\left(\begin{array}{cc}0&1\\ 0&0\\\end{array}   \right)$, then there is
no K-contact structure on {\Goth g}.

\section{Complex contact structures}

This section deals solely with complex contact Lie groups, i.e., complex Lie groups with a left-invariant holomorphic 1-form $\eta$ such that $\eta \wedge d\eta^n\neq 0,$ where the complex dimension of the Lie group is $2n+1.$  Within this section, we will use the same notation for the resulting structures and forms in the complex contact Lie theoretical category as we did in the real category.  So, like the real case, we let $\HH$ be left-invariant distribution given as the kernel of $\eta$ (in the holomorphic tangent bundle) and $\xi$ be the left-invariant vector field given by $\eta(\xi)=1$ and $d\eta(\xi,*)=0.$  It is only in the next section, where we are using both real and complex contact structures simultaneously that we will use different notation for the different categories.  This material has already been published across two papers, \cite{foreman1} and \cite{foreman2}, but for completeness and coherence we provide here a unified and streamlined presentation of the relevant results.

Suppose $(G,\eta)$ is a $(2n+1)$-dimensional complex contact Lie group
such that the adjoint representation of the Reeb vector field $\xi$ on the Lie algebra of $G$, $\mgoth{g}$
is diagonalizable.
Let $A\subset \C$ be the set of all
eigenvalues of $ad(\xi)$ with nontrivial eigenvectors.  We call $A$ the {\it roots of $\xi$}.  For each
$\alpha\in \C,$ set
$$ \mgoth{g}_\gra=\{X\in \mgoth{g}: [\xi,X]=\gra X\ \} $$
For $X\in \mgoth{g},$ $X=c\xi + \HH X$ for some $c\in\C$ so that
$[\xi,X]=[\xi,c\xi + \HH X]=[\xi,\HH X].$  Thus, if $\gra\in
A-\{0\},$ then $\mgoth{g}_\gra\subset \HH.$

\begin{proposition}\label{fundamental prop}
Let $(G,\eta)$ be a $(2n+1)$-dimensional complex contact Lie group
such that the adjoint representation of the Reeb vector field $\xi$
is diagonalizable with roots given by the set $A$. Then 
\begin{enumerate}
\item If $X\in\mgoth{g}_\gra$ and
$Y\in\mgoth{g}_\beta$ for $\gra,\ \beta\in A,$ then $ad(\xi)[X,Y]=(\gra+\beta)[X,Y]$ and
either $\gra+\beta=0$ or $d\eta(X,Y)=0.$

\item For any $\gra\in A$ and $X\in
\mgoth{g}_\gra-(0),$ there is a $Y\in \mgoth{g}_{-\gra}$ such that
$[X,Y]=\xi +Z$ for some $Z\in \mgoth{g}_0\cap\HH.$
\end{enumerate}
\end{proposition}
\vskip .25in

\noindent{\it Proof:}
For Statement 1, the Jacobi identity gives us:
\begin{eqnarray*}
0&=&[[\xi,X],Y]+[[X,Y],\xi]+[[Y,\xi],X]\\
&=&\gra[X,Y]+[[X,Y],\xi]-\beta[Y,X].\\
\end{eqnarray*}
So, $ad(\xi)[X,Y]=(\gra+\beta)[X,Y].$
In particular,\ $\eta(ad(\xi)[X,Y])=(\gra+\beta)\eta([X,Y]).$ By
definition of $\xi,$ the left-hand side is zero.  Furthermore,
$\eta([X,Y])=-2d\eta(X,Y).$  This proves Statement 1.

Let $\gra\in A$ and
$X\in \mgoth{g}_\gra-(0).$  Since $d\eta^n\neq0$ on $\HH$, we know
that there exists $Y\in \HH$ such that $[X,Y]=\xi+Z$ for some $Z\in
\HH.$ In fact, if we create a basis of $\HH$ such that each element
of the basis is an eigenvector of $ad(\xi),$ we see that there is
some $\beta\in A$ such that $Y\in \mgoth{g}_{\beta}$ and
$[X,Y]=\xi+Z$ for some $Z\in \HH.$ By Statement 1, $\beta=-\gra.$  Also,
$0=ad(\xi)([X,Y])=ad(\xi)(Z).$  This proves Statement 2.

\vskip .25in

\begin{theorem}\label{adxi=0}
Let $(G,\eta)$ be a $(2n+1)$-dimensional complex contact Lie group
such that the adjoint representation of the Reeb vector field $\xi$
is diagonalizable.  If $n>1,$ then $ad(\xi)=0.$

\end{theorem}

\vskip .25in \noindent {\bf Proof:}\ \ We prove this theorem by systematically reviewing the cases where $A\neq\{0\}$ and showing that each such possible case creates a contradiction.  First, we consider the situation in which $ad(\xi)$ has no zero eigenvectors in $\HH$ and two distinct nonzero eigenvalues, $\gra$ and $\beta\neq -\gra.$  Second, we investigate the case in which $ad(\xi)$ has exactly two eigenvectors in $\HH$, $\gra\neq 0$ and $-\gra.$  Finally, we consider the situation in which both $\gra\neq 0$ and $0$ are eigenvalues of $ad(\xi)$ in $\HH.$  We will show that each of these cases lead to a contradiction.

\noindent{\it Case 1:}\ \  Assume that $ad(\xi)$ has no zero eigenvectors in $\HH$ and two distinct nonzero eigenvalues, $\gra$ and $\beta\neq -\gra.$  Without losing any generality, we can assume that $\gra\pm\beta\notin A.$  In particular, by Proposition \ref{fundamental prop}, $-\gra\in A$, and $\left[X_\gra, \mgoth{g}_{-\gra}\right]=\left<\xi\right>$ for any $X_\gra\in \mgoth{g}_\gra.$  Furthermore, $[\mgoth{g}_{\pm\gra},\mgoth{g}_\beta]=(0).$   

Let $X_\gra\in \mgoth{g}_\gra,\ X_\beta\in\mgoth{g}_\beta,$ both non-zero.  By the Jacobi identity,
\begin{eqnarray*}
\beta X_\beta&=&ad(\xi)X_\beta\\
&=&ad\left([X_{\gra},X_{-\gra}]\right)X_\beta\\
&=&(ad X_\gra)(ad X_{-\gra})(X_\beta)-(ad X_{-\gra})(ad X_{\gra})(X_\beta)\\
&=&(ad X_\gra)([ X_{-\gra},X_\beta]-(ad X_{-\gra})[X_{\gra},X_\beta]\\
&=&0,
\end{eqnarray*}
since $[\mgoth{g}_{\pm \gra}, \mgoth{g}_\beta]=(0).$  Thus, $\beta=0,$ a contradiction.
\vskip .25in
\noindent{\it Case 2:}\ \ Assume that $ad(\xi)$ has exactly two eigenvectors in $\HH$, $\gra\neq 0$ and $-\gra.$  Let $\underline{E}=\{E_1,\dots, E_{2n}\}$ be a basis of $\HH$ such that
\begin{eqnarray*}
\mgoth{g}_{-\gra}&=&\left<E_{2j-1}:j=1,\dots, n\right>\\
\mgoth{g}_\gra&=&\left<E_{2j}:j=1,\dots, n\right>,\\
\end{eqnarray*}
that is, $ad(\xi)E_k=(-1)^k\gra E_k$ for $k=1, \dots, 2n.$
By Proposition \ref{fundamental prop}, $ad(\xi)[E_k, E_l]=\left((-1)^k+(-1)^l \right)\gra[E_k, E_l].$  In particular, $ad(\xi)[E_{2j_1},E_{2j_2-1}]=0$ and $0=[E_{2j_1},E_{2j_2}]=[E_{2j_1-1},E_{2j_2-1}]$ for $j_1,\ j_2=1,\dots, n$  (since $\mgoth{g}_{\pm 2\gra}=(0)$ by assumption).  Thus, since $\mgoth{g}_0=<\xi>$ by assumption, for each $k,\ l=1, \dots, 2n,$ $[E_k, E_l]=\beta_{kl} \xi$ for some $\beta_{kl}\in \C$ with $0=\beta_{\rm even \ even}=\beta_{\rm odd\ odd}.$  Furthermore, the fact that $\HH$ is a complex contact structure on $G$ implies that for every $k=1, \dots, 2n$, there is a $\tilde k=1, \dots, 2n$ such that $\beta_{k\tilde k}\neq 0.$  Without loss of generality, we can assume that $\beta_{2j2j-1}\neq 0$ for $j=1, \dots, n.$

Then
\begin{eqnarray*}
0&=&[[E_1,E_2],E_3]+[[E_2,E_3],E_1]+[[E_3,E_1],E_2]\\
&=&\beta_{12}[\xi,E_3]+\beta_{23}[\xi,E_1]\\
&=&-\gra\beta_{12} E_3-\gra \beta_{23} E_1.\\
\end{eqnarray*}
Thus, $\gra=0,$ which contradicts the assumption that $\gra\neq 0.$
\vskip .25in
\noindent{\it Case 3:}  Assume that both $\gra\neq 0$ and $0$ are eigenvalues of $ad(\xi)$ in $\HH.$

Proposition \ref{fundamental prop} implies that $[\mgoth{g}_0,\mgoth{g}_0]\subset
\mgoth{g}_0.$  Let $X_1$ be a nonzero element of $\mgoth{g}_0\cap \HH$.
Then, again by Proposition \ref{fundamental prop}, there is an element $\tilde X_2\in \mgoth{g}_0\cap\HH$ such that
$\eta\left([X_1,\tilde X_2]\right)\neq 0.$ By considering the Jordan canonical form of $ad(X_1)$ restricted on $\mgoth{g}_0 $, we see that there is an $X_2\in \mgoth{g}_0\cap \HH$ such
that $[X_1,X_2]=\xi.$ Furthermore, 
$ad(X_j)(\mgoth{g}_\gra)\subset \mgoth{g}_\gra$ for each $j=1,2.$
The Jacobi identity implies that
$\left[ad(X_1),ad(X_2)\right]=ad([X_1,X_2])=ad(\xi)$ so that, on
$\mgoth{g}_\gra,$ $\left[ad(X_1),ad(X_2)\right]=\gra I.$  But, for
any linear transformations $S$ and $T$ on a given vector space $V,$
$ST-TS$ is never a non-zero multiple of the identity.  Thus, we have
a contradiction.  Having exhausted all possibilities in which $A\neq\{0\},$ we have proven the theorem.

\section{Main theorem}

We now prove the main result as an easy corollary of Theorem \ref{adxi=0}.  

\begin{theorem} Any K-contact Lie group of dimension five or greater
is the central extension of a symplectic Lie group. \end{theorem} 

{\bf
Proof:}\ \ Given a real contact Lie algebra $(\mgoth{g},\eta)$, the
complexification $\mgoth{g}^\C$ is a complex contact Lie algebra with complex contact form given by $\eta^\C(X+iY)=\eta(X)+i\eta(Y)$ for $X,\ Y\in \mgoth{g}.$  The complex Reeb vector field $\xi^\C$ in $\mgoth{g}^\C$ is defined by:
$$\eta^\C(\xi^\C),\ \ d\eta^\C(\xi^\C,\ *)=0.$$
Since $\xi\in \mgoth{g}\subset \mgoth{g}^\C$ satisfies this condition, $\xi^\C=\xi.$
Thus, the adjoint operator $ad(\xi^\C)$ is simply the complex extension of
$ad(\xi)$ on $\mgoth{g}$ acting on $\mgoth{g}^\C$. 

In addition, suppose that $g$ is a left-invariant associated metric on $\mgoth{g}$ such that 
$(\mgoth{g},\eta, \xi, g)$ is a K-contact Lie algebra.  There is then an orthonormal basis $\underline{e}$ of 
$\mgoth{g}$ with respect to which the matrix representation of $ad(\xi)$ is skew-symmetric.  

Then the operator $ad(\xi^\C)$ is diagonalizable on $\mgoth{g}^\C$ with purely imaginary eigenvalues.  By
Theorem \ref{adxi=0}, $ad(\xi^\C)=0$, which implies that $ad(\xi)=0.$  And so $(\mgoth{g}, \eta, \xi)$ is the
central extension of a symplectic Lie algebra.

\end{document}